\documentclass[letterpaper,reqno]{amsart}

\usepackage{amsmath} 
\usepackage{amsfonts} 
\usepackage{amssymb} 
\usepackage{epsfig} 
\usepackage{graphicx} 
\usepackage{multirow} 
\usepackage{verbatim} 
\usepackage{rotating} 
\usepackage[all]{xy} 
\usepackage{appendix}

\newtheorem{theorem}{Theorem}[section]

\theoremstyle{definition}
\newtheorem{definition}[theorem]{Definition}

\theoremstyle{remark}
\newtheorem{remark}[theorem]{Remark}
\theoremstyle{notation}

\numberwithin{equation}{section}
\theoremstyle{corollary}

\newtheorem{conjecture}[theorem]{Conjecture}
\usepackage[bookmarks=false]{hyperref}
\newcommand{\Map}{\mathrm{Map}}

\newcommand{\Cat}{\mathsf{Cat}}

\newcommand{\coh}{\mathrm{H}}

\begin{document}

\title[]{Analogy between the cyclotomic trace map $K\rightarrow TC$ and the Grothendieck trace formula via Noncommutative Geometry}

\author[]{Ilias Amrani}
\address{Academic University of St-Petersburg, Russian Federation.}
\email{ilias.amranifedotov@gmail.com}


\subjclass[2000]{ 19E08, 14F30, 14A22, 55N15, secondary 11M38}

\dedicatory{Dedicated to my Father}

\keywords{Noncommutative Algebraic Geometry, Grothendieck trace formula, Cyclotomic trace map, Algebraic $K$-Theory, Topological Cyclic Homology, Crystalline Cohomology, l-adic Cohomology, Categorification, Rational points, Homotopy fixed points. }

\begin{abstract}
In this article, we suggest a categorification procedure in order to capture an analogy between Crystalline Grothendieck-Lefschetz trace formula and the cyclotomic trace map $K\rightarrow TC$ from the algebraic $K$-theory to the topological cyclic homology $TC$. First, we categorify the category of schemes to the $(2,\infty)$-category of noncommuatative schemes \`a la Kontsevich. This gives a categorification of the set of rational points of a scheme. Then, we categorify the Crystalline Grothendieck-Lefschetz trace formula and find an analogue to the Crystalline cohomology in the setting of noncommuative schemes over $\mathbf{F}_{p}$. Our analogy suggests the existence of a categorification of the $l$-adic cohomology trace formula in the noncommutative setting for $l\neq p$. Finally, we write down the corresponding dictionary.   
\end{abstract}
\maketitle
\section{Arithmetical side: Grothendieck trace formula}

In this short expository article, we explore the notion of \textit{noncommutative algebraic space} and use a \textit{categorification and stabilization procedures} to make a surprising analogy between two formulas arising from, a priori, different areas of mathematics.
\begin{enumerate}
\item Crystalline and $l$-adic Gorthendieck-Lefschetz trace formula for smooth $\mathbf{F}_{p}$-schemes.
\item The cyclotomic trace map $K\rightarrow TC$ from the algebraic $K$-theory to the topological cyclic homology for smooth $\mathbf{F}_{p}$-schemes.
\item The bridge between the two precedent trace formulas is given via the categorification procedure from the category of schemes to the $(2,\infty)$-category of dg-categories (Morita equivalences), sending any scheme $X$ to the dg-category of perfect complexes $\mathsf{Perf}(X)$. 
\end{enumerate}
\subsection*{Grothendieck-Lefschetz trace formula: $l$-adic cohomology}
For more details we refer to \cite{milne1980etale}. We fix a finite field $\mathbf{F}_{p}$ with $p$ a prime and let $X$ be a \textit{nice} scheme. An important challenge in arithmetical algebraic geometry is counting the number of rational points of $X$, namely the set $X(\mathbf{F}_{p})$. Grothendieck has defined an appropriate cohomology theory ($l$-adic cohomology with compact support, $l$ different from $p$) in order to mimic the famous Lefschetz fixed point 
\begin{eqnarray}\label{grothc}
\sharp\{X(\mathbf{F}_{p^{n}})\} = \sum_{i=0}(-1)^{i}\mathbf{trace}[F^{n}: \mathrm{H}_{\textit{\'et}, c}^{i}(\overline{X}, \mathbf{Q}_{l})], l\neq p.
\end{eqnarray} 
where $\overline{X}=X\otimes_{ \mathbf{F}_{p}}\overline{\mathbf{F}}_{p}$ and  $F^{n}: \overline{X}\rightarrow \overline{X}$ is the $n^{th}$-iterated Frobenius morphisms. 
In particular when $X$ is a smooth regular scheme of dimension $d$ the formula simplifies
\begin{eqnarray}\label{grothsm}
\sharp \{X(\mathbf{F}_{p^{n}})\} = \sum_{i=0}^{2d}(-1)^{i}\mathbf{trace}[F^{n}: \mathrm{H}_{\textit{\'et}}^{i}(\overline{X}, \mathbf{Q}_{l})], l\neq p. 
\end{eqnarray}
The cohomology $\coh_{\textit{\'et}}^{\ast}(\overline{X},\mathbf{Q}_{l})$ is defined by the formula 
\begin{eqnarray}\label{ladic}
\coh_{\textit{\'et}}^{\ast}(\overline{X},\mathbf{Q}_{l}):=\lim_{k}\coh_{\textit{\'et}}^{\ast}(\overline{X}, \mathbf{Z}/l^{k}\mathbf{Z}) \otimes \mathbf{Q}_{l}.
\end{eqnarray}
We denote the $l$-adic cohomology 
\begin{eqnarray}\label{2}
\coh_{l-adic}^{\ast}(X) :=\lim_{k}\coh_{\textit{\'et}}^{\ast}(\overline{X}, \mathbf{Z}/l^{k}\mathbf{Z})
\end{eqnarray}
and recall that 
\begin{eqnarray}\label{fixfrob}
\sharp \{X(\mathbf{F}_{p^{n}})\} = \sharp \{x\in X(\overline{\mathbf{F}}_{p})| F^{n}(x)=x\}.
\end{eqnarray}
It seems that the origin of the \'etale site and $l$-adic cohomology is closely related to realizing the analogue of Lefschetz formula for $\mathbf{F}_{p}$-schemes. However, we notice that Grothendieck formula has a limitation namely the restriction $p\neq l$. The case $p=l$ seems to suggest the existence of an other appropriated cohomology theory. 
   
One of our main goal is to propose an analogue of $l$-adic cohomology (and Crystalline cohomology) in the context of algebraic noncommutative spaces which can be interpreted as \textit{noncommutative intersection theory}. 
\subsection*{Grothendieck-Lefschetz trace formula: Crystalline cohomology}
For more details, we refer the reader to \cite{katz1981crystalline} and \cite{chambert1998cohomologie}. As we noticed, the Grothendieck-Lefschetz trace formula for (enough good) $\mathbf{F}_{p}$-scheme $X$ works when we use the $l$-adic cohomology for $l\neq p$. The natural question to ask is the following: What happens if $l=p$ ? do we still have a trace formula?  The answer to this question is the Crystalline cohomology theory $\coh_{cris}(-)$. 
\begin{eqnarray}\label{cris}
\sharp\{X(\mathbf{F}_{p^{n}})\}= \sum_{i=0}^{\infty}(-1)^{i}\mathbf{Trace}[F^{n}: \coh_{cris}(X)\otimes_{\mathbf{Z}_{p}} \mathbf{Q}_{p}] 
\end{eqnarray}
where $\mathbf{Q}_{p}$ is fraction field of the total ring of Witt vectors $\lim_{k}W_{k}(\mathbf{F}_{p})\simeq W(\mathbf{F}_{p})\simeq \mathbf{Z}_{p}$ ($p$-adic integers)
and 
\begin{eqnarray}\label{cris1}
\coh^{i}_{cris}(X)= \lim_{k}\coh^{i}_{cris}(X/W_{k}(\mathbf{F}_{p}))
\end{eqnarray}
\subsection*{Relation to de Rham-Witt complex}
The relation between de Rham-Witt cochain complex \cite{illusie1979complexe} and Crystalline cohomology is explained in \cite[section 3 and 4]{chambert1998cohomologie}, the hypercohomolgy of de Rham-Witt complex $W\Omega^{\ast}(-)$ verifies the following property for (good enough) $\mathbf{F}_{p}$-scheme $X$, namely the isomorphism  
\begin{eqnarray}\label{derhamwitt}
  \coh^{i}_{cris}(X)\cong \mathbb{H}^{i}(X, W\Omega^{\ast}(-))
  \end{eqnarray}  

\section{Homotopical side: Cyclotomic trace formula $K\rightarrow TC$}

 The other side of the picture comes from a homotopical approximation of the Algebraic $K$-theory. For more details and comprehension we refer to \cite{dundas2012local}. It is well known that the algebraic $K$-theory of schemes gives very important arithmetical informations. For any scheme $X$ we associate the the differential graded category of perfect complexes $\mathsf{Perf}(X)$, we have a trace map from the spectrum of algebraic $K$-theory to the topological Hochschild homology
 \begin{eqnarray}\label{tr}
K(\mathsf{Perf}(X)):=K(X)\rightarrow THH(X):=THH(\mathsf{Perf}(X)).
\end{eqnarray}
The $THH$ has a natural action of the topological circle $S^{1}$ and the trace map factors trough the homotopy fixed points 
\begin{eqnarray}\label{trcnaive}
K(X)\rightarrow THH(X)^{S^1}\rightarrow THH(X)^{hS^{1}}\rightarrow THH(X).
\end{eqnarray}
The problem is that the construction $THH(X)^{S^1}$ is not homotopy invariant and $THH(X)^{hS^{1}}$ is not a very good approximation in general. There is an other alternative which consists to construct an other spectra $TH(-)$ with the following properties \cite{May}:
\begin{itemize}
\item  $TH(-)$ is equivalent to $THH(-)$ as a naive $S^{1}$-spectra.
\item  It is a cyclotomic $S^{1}$-spectra. 
\item  It comes equipped with restriction map $R$, Frobenius map $F$ and Verschiebung map $V$ (wich we don't use explicitly).
\end{itemize} 
Recall that rationally $THH(X)^{hS^{1}}$ coincides with \textit{Cyclic Homology}. Now, we fix a prime $l$ and we take an infinite increasing chain of inclusions of finite cyclic groups 
$$ C_{l}\subset C_{l^{2}}\subset C_{l^{3}}\subset \dots S^{1}$$ and define a new homology theory \cite{blumberg2012localization}
\begin{eqnarray}\label{tr}
TR(X,l)=holim_{n}TH(X)^{C_{l^{n}}}
\end{eqnarray}
using the restriction map $R:TH(X)^{C_{l^{n+1}}}\rightarrow TH(X)^{C_{l^{n}}}$. If the natural map $TH(X)^{C_{l}}\rightarrow TH(X)^{hC_{l}}$ is an equivalence after $l$-completion, then by Tsalidis theorem  $TH(X)^{C_{l^{n}}}\rightarrow TH(X)^{hC_{l^{n}}}$ is also an equivalence after $l$-completion for any $n>0$. This argument was intensively used in \cite{rognes1999topological} in order to compute the 2-primary part of the $K$-theory of integers. In this case $TR(\mathbf{Z},2)$ is equivalent to $THH(\mathbf{Z})^{hS^{1}}$ after $2$-completion.  
Now, we are ready to give the definition of the topological cyclic homology of $X$ at prime $l$.
\begin{definition}
The topological cyclic homology of $X$ at prime $l$ is given by \cite{hesselholt1993stable}
\begin{eqnarray}\label{tc}
TC(\mathsf{Perf}(X),l)=TC(X,l):= TR(X,l)^{hF}= TR(\mathsf{Perf}(X),l)^{hF},
\end{eqnarray}
where $F:TR(X,l)\rightarrow TR(X,l)$ is the Frobenius map and $TR(X,l)^{hF}$ is the spectra of the homotopy fixed points with respect to $F$.  
\end{definition}
\begin{remark}
The derived spectrum of maps from the sphere spectrum $S$ to $TH(X)$ in the homotopy category of $l$-cyclotomic spectra is naturally isomorphic to $TC(X,l)$ after $l$-completion \cite[Theorem 1.4]{blumberg2013homotopy}.
\end{remark}
\begin{definition}
The trace map $K(X)\rightarrow TH(X)$ is $S^{1}$-equivariant (the $S^1$-action on $K(X)$ is trivial), we have a natural map 
\begin{eqnarray}
trc: K(X)\rightarrow TC(X,l):=TR(X,l)^{hF}
\end{eqnarray}
called the cyclotomic trace map. 
\end{definition}
When $p=l$, and $X$ is an affine smooth scheme over $\mathbf{F}_{p}$, it turns out that $trc$ is a very good approximation (after $p$-completion of the connective covers).

\subsection{ Categorification of Euler characteristic and Lefschetz formula}

Let $X$ be any finite CW complex $X$ of dimension $n$, the Euler characteristic is defined as 
\begin{eqnarray}\label{euler}
\chi(X)= \sum_{i=0}^{n} (-1)^{i} dim_{\mathbf{Q}}H^{i}(X,\mathbf{Q}).
\end{eqnarray}
It follows that the Euler characteristic is defined only by using the numerical invariant i.e. the dimension of the vector spaces of the rational cohomolgy $H^{i}(X,\mathbf{Q})$. A generalization of the Euler characteristic is the Lefschetz formula (in the case where $X$ is a triangulated manifold of dimension $n$) given by 
\begin{eqnarray}\label{lef}
\Lambda(X,F):= \sum_{i=0}^{n} (-1)^{i} \mathbf{trace}[F^{\ast}:H^{i}(X,\mathbf{Q})],
\end{eqnarray}
where $F:X\rightarrow X$ is continuous endomorphism of $X$. If $\Lambda(X,F)\neq 0$ then $F$ has a least one fixed point. More precisely when the set of fixed point $Fix(F)$ is finite then
$$ \Lambda(X,F)= \sum_{x\in Fix(F)} i(F,x) $$
where $i(F,x)$ is the index (integer) of $F$ at $x$. 

We should remark that  $\Lambda(X,Id)=\chi(X)$. A reasonable categorification of Euler formula is the singular cochain complex $C^{\ast}(X,\mathbf{Q})$ or more artificially (at least a priori) we can say that the categorification of Euler characteristic is $C^{\ast}(X,\mathbf{Q})\simeq C^{\ast}(X,\mathbf{Q})^{h Id}$, this suggests that a natural categorification of the right side of the Lefschetz formula is given by 
\begin{eqnarray}\label{catlef}
C^{\ast}(X,\mathbf{Q})^{h F},
\end{eqnarray}
and the decategorification is given by 
\begin{eqnarray}
\sum_{i=0}^{n}\mathbf{trace}[F^{\ast}: \coh^{i}(X,\mathbf{Q})].
\end{eqnarray}

\section{The Bridge: Noncommutative algebraic geometry}
The previous categorification seems to suggests that the Grothendieck-Lefschetz formula should be categorified in the following way.   
For any scheme $X$ over $\mathbf{F}_{p}$, the set of rational points $X(\mathbf{F}_{p})$ is corepresentable in the category of $\mathbf{F}_{p}$-schemes, indeed
\begin{eqnarray}\label{rationalpoints}
\Map_{Sch}(spec(\mathbf {F}_{p}), X)= X(\mathbf {F}_{p}).
\end{eqnarray}
Kontsevich has proposed a noncommutative point of view of schemes which was formalized later by Tabuada. Let $\mathsf{dgCat}_{k}$ denotes the model category of small dg-categories and dg-functors where the weak equivalences are Morita equivalences \cite{tabuada2011guided}. For any scheme $X$ can associate the differential graded category $\mathsf{Perf}(X)$ which is the full dg-subcategory of compact object in enhanced derived category $D^{dg}(X)$. The model category $\mathsf{dgCat}_{k}$ is actually a symmetric monoidal $(2,\infty)$-category. Notice that the existence of the derived internal $\mathsf{RHom} $ in $\mathsf{dgCat}_{k}$ is due to To\"en \cite{toen2007homotopy}. Any differential graded $k$-algebra $A$ can be viewed as a differential graded category with one object in obvious way and its corresponding fibrant replacement is the dg-category of perfect complexes $\mathsf{Perf}(A)$. According to Tabuada, the model category $\mathsf{dgCat}_{k}$ \textbf{can be stabilized  and localized} such that we have a universal stabilization functor 
$$U_{add}: \mathsf{dgCat}_{k}\rightarrow \mathsf{Mot}_{add}$$ where 
$\mathsf{Mot}_{add}$ is stable model category (of additive noncommutative motives) and for any dg-categories $\mathsf{A}$ and $\mathsf{B}$, the derived spectral enrichment 
$$\Map_{\mathsf{Mot}_{add}}(U_{add}(\mathsf{A}), U_{add}(\mathsf{B}))$$ 
is isomorphic to Whaldhausen $K$-theory 
$$ K(\mathsf{RHom}(\mathsf{A},\mathsf{B}))$$
 in the homotopy category of spectra. In particular when $\mathsf{A}\simeq \mathsf{Perf}(k)$ and \\
 $\mathsf{B}=\mathsf{Pref}(X) $ then 
\begin{eqnarray}\label{motivicpoints}
\Map_{\mathsf{Mot}_{add}}(U_{add}(\mathsf{Perf}(k)), U_{add}(\mathsf{Perf}(X)))\simeq K(\mathsf{Perf}(X)):=K(X).
\end{eqnarray}
Therefore, the algebraic $K$-theory is corepresentable by $\mathsf{Perf}(k)$ \cite[Theorem 6.1]{tabuada2011guided}. It is natural to suggest that the categorification of $X(\mathbf{F}_{p})$ \ref{rationalpoints} is given by the ($E_{\infty}$-ring) spectra $K(\mathsf{Pref}(X)):=K(X)$ \ref{motivicpoints}.

Now, the idea is to find the the analogue of the $l$-adic cohomology of schemes in the context of the noncommutative geometry. First, we should notice that the $l$-adic Grothendieck trace formula is valid when $l\neq p$ while in the noncommutative side of the picture the cyclotomic trace map is  a very good approximation when $l=p$. It becomes clear ones we understand the meaning of the homotopy groups of $TR(X,p)^{\wedge}_{p}$ for a smooth $\mathbf{F}_{p}$-scheme $X$. It was proved \cite[page 18]{geisser1999topological} that there is an isomorphism of cochain complexes 
\begin{eqnarray}\label{criswitt}
 \pi_{\ast}TR(-,p)^{\wedge}_{p}\cong W\Omega^{\ast}(-).
\end{eqnarray}
For the cochain complex structure on $\pi_{\ast}TR(X,p)^{\wedge}_{p}$ we refer to\cite{May}. We conclude that the cohomology theory $TR(-,p)^{\wedge}_{p}$ is the presheaf on the category of smooth $\mathbf{F}_{p}$-schemes of de Rham-Witt cochain complex. The decategorification of the algebraic $K$-theory $K(X)$ is $X(\mathbf{F}_{p})$ and the decategorification of  $TR(X,p)^{ \wedge hF}_{p}:=TC(X,p)^{\wedge}_{p}$  gives the right side of Crystalline Grothendieck trace formula by applying  the hypercohomology i.e.
\begin{eqnarray}\label{hyp}  
\sum_{i=0}(-1)^{i}\mathbf{Trace}[F:\mathbb{H}^{i}(X, \pi_{\ast}TR(-,p)^{\wedge}_{p})\otimes_{\mathbf{Z}_{p}} \mathbf{Q}_{p}]
\end{eqnarray}
which is equivalent to
$$\sum_{i=0}(-1)^{i}\mathbf{Trace}[F:\coh^{i}_{cris}(X)\otimes_{\mathbf{Z}_{p}} \mathbf{Q}_{p}],$$
since $ \pi_{\ast}TR(-,p)^{\wedge}_{p}\cong W\Omega^{\ast}(-)$ \ref{criswitt} and $  \coh^{i}_{cris}(X)\cong \mathbb{H}^{i}(X, W\Omega^{\ast}(-))$ \ref{derhamwitt}.
 Hence, by decategorification of the cyclotomic trace map $K(X)\rightarrow TR(X,p)^{\wedge hF}_{p}=TC(X,p)^{\wedge}_{p}$ we recover the Crystalline Grothendieck trace formula \ref{cris}!  
\section{Analogy and Conjectures}\label{analogy}
The passage from the commutative world (schemes) to the noncommutative one (dg-categories) is given by the functor
$$X\mapsto \mathsf{Perf}(X)$$
followed by the universal additive stabilization (in the sense of Tabuada) 
$$  \mathsf{Perf}(X)\mapsto U_{add}\mathsf{Perf}(X).$$
\begin{tabular}{||l|c||}
\hline
 category of smooth schemes over $\mathbf{F}_{p}$ & $\infty$-stable category $\mathsf{Mot}_{add}$ over $\mathbf{F}_{p}$ \\
 \hline \hline 
 $X$ & $U_{add}\mathsf{Perf}(X)$ \\
 \hline
 Rational points & Motivic points \\
 $\Map_{Sch}(spec(\mathbf {F}_{p}), X)= X(\mathbf {F}_{p})$ &
$Map_{\mathsf{Mot}_{add}}(U_{add}(\mathsf{Perf}(\mathbf{F}_{p})), U_{add}(\mathsf{Perf}(X)))\simeq K(X) $ \\

 & \\
\hline
Crystalline cohomology $H^{\ast}_{cris}(X) $ & cohomology theory $ TR(X,p)^{\wedge}_{p}$\\

\hline
& \\

$l$-adic cohomology  $H^{\ast}_{l-adic}(X)$ &   conjecture \ref{conjecture}\\
\hline
&\\
$ \sum_{i=0}(-1)^{i}\mathbf{trace}[F: \mathrm{H}_{cris}^{i}(X)\otimes\mathbf{Q}_{p}]$  & $ TR(X,p)^{\wedge hF}_{p}:=TC(X,p)^{\wedge}_{p}$ \\
              
\hline 
Crystalline Grothendieck trace formula &  cyclotomic trace map \\
$\sharp \{X(\mathbf{F}_{p})\}=\sum_{i=0}(-1)^{i}\mathbf{trace}[F: \mathrm{H}_{cris}^{i}(X)\otimes\mathbf{Q}_{p}]$ &  $K(X)^{\wedge}_{p}\rightarrow TC(X,p)^{\wedge}_{p}$\\
 right side = left side  & cyclotomic trace map is an equivalence, conjecture \ref{prediction}\\  
 \hline
\end{tabular}
\\\\
We still don't have an analogue of the $l$-adic cohomology for noncommutative spaces, it was conjectured by Quillen and Lichtenbaum that the canonical map from  $K$-theory to the \'etale $K$-theory is a very good approximation for (good enough) scheme $X$. More precisely, the natural map of spectra 
 $$ K(X)\rightarrow K^{\textit{\'et}}(X)$$
 induces an isomorphism on (stable) homotopy groups $\pi_{i}$ for sufficiently large $i$ (related to the dimension of $X$). From our previous discussion, we can say that this approximation does not fit in the spirit of our analogy. We propose the following alternative conjecture.
 \begin{conjecture}\label{conjecture}
There exists a (co)homology theory $\mathbb{TR}$ for dg-categories with the following properties:
\begin{enumerate}
\item the functor $\mathbb{TR}: \mathsf{dg}\Cat_{\mathbf{F}_{p}}\rightarrow \mathsf{Sp}$ is Morita invariant i.e. if $\mathsf{C}\rightarrow \mathsf{D}$ is a Morita equivalence in $\mathsf{dg}\Cat_{\mathbf{F}_{p}}$, then $\mathbb{TR}(\mathsf{C})\rightarrow \mathbb{TR}(\mathsf{D})$ is a stable equivalence of spectra.
\item The (co)homology theroy $\mathbb{TR}(-)$ is a localizing invariant in the sense of Tabuada \cite[definition 5.1]{tabuada2011guided}.  
\item There is a natural transformation $K(-)\rightarrow \mathbb{TR}(-)$.
\item For any $\mathsf{C}\in \mathsf{dg}\Cat_{\mathbf{F}_{p}}$, then $\pi_{\ast}\mathbb{TR}(\mathsf{C})$ has a natural structure of cochain complex.
\item If $X$ is (good enough) $\mathbf{F}_{p}$-scheme and a prime $l$ different form $p$ then :
\begin{itemize} 
\item the hypercohomology $\mathbb{H}^{i}(X,\pi_{\ast}\mathbb{TR}(\mathsf{Perf}(-))^{\wedge}_{l})$ is isomorphic to the $l$-adic cohomology $\coh_{l-adic}^{i}(X)$.
\item  There is a Frobenius natural transformation $F: \mathbb{TR}(\mathsf{Perf}(-))\rightarrow \mathbb{TR}(\mathsf{Perf}(-))$
\item It recovers the $l$-adic trace formula a la Lefschetz-Grothendieck.
\item The natural transformation $K(\mathsf{Perf}(-))\rightarrow \mathbb{TR}(\mathsf{Perf}(-))$ factors as trace map $trc: K(\mathsf{Perf}(-))\rightarrow \mathbb{TR}(\mathsf{Perf}(-))^{hF}$, and this new trace map is a good (connective) approximation after $l$-completion for $l\neq p$. More precisely, 
$$\pi_{i}K(\mathsf{Perf}(X))\otimes \mathbf{Z}_{l}\cong \pi_{i}(\mathbb{TR}(\mathsf{Perf}(X))^{hF}) \otimes \mathbf{Z}_{l}$$
for $i\geq 0$.
\end{itemize}

\end{enumerate}
\end{conjecture} 
Since the Crystalline trace formula \ref{cris} is true for any nice $\mathbf{F}_{p}$-scheme $X$ the analogy has a very natural prediction which we formulate in the following conjecture:
\begin{conjecture}\label{prediction}
For any nice $\mathbf{F}_{p}$-scheme $X$, the cyclotopmic trace map $K(X)\rightarrow TC(X,p)$ is a weak equivalence (between connective covers) after $p$-completion.
\end{conjecture}
\section{Secondary K-theory and $\mathbf{F}_{p^{k}}$-points}
Based on our analogy \ref{analogy}, we want to illustrate the idea that the algebraic $K$-theory of a nice $\mathbf{F}_{p}$-scheme $X$ contains a substantial informations about the concrete set of rational points $X(\mathbf{F}_{p})$. In order to establish the categorification procedure from rational points $X(\mathbf{F}_{p})$ to $K(X)$, we have used the fact that the point functor is corepresentable in the category of schemes by $ spec (\mathbf{F}_{p})$ and the algebraic $K$-theory is corepresentable by the noncommutative motive $U_{add}\mathsf{Perf}(\mathbf{F}_{p})$.
\begin{definition}\label{secondary}
The $k$-secondary algebraic $K$-theory for $\mathbf{F}_{p}$-schemes is defined as $$K^{(k)}(X)=TR(X,p)^{ hF^{k}}$$ where $F^{k}$ is the $k$-iterated Frobenius operator.
\end{definition}
 A natural question pop-up:
$$\textbf{Can we recover $X(\mathbf{F}_{p^k})$ form $K^{(k)}(X)$ ? } $$ 

If $k=1$, a short answer is \textit{yes} if the conjecture \ref{prediction} is true by using \ref{hyp} and the Crystalline Grothendieck-Lefschetz formula \ref{cris}. It seems that we can recover $X(\mathbf{F}_{p^{k}})$ using the decategorification of the homotopy fixed points $TR(X,p)^{\wedge hF^{k}}_{p}$ with respect to the $k$-iterated Frobenius operator $F^{k}$ i.e.,

$$\sum_{i=0}(-1)^{i}[F^{k}:\mathbb{H}^{i}(X,\pi_{\ast}TR(-,p)^{\wedge}_{p}\otimes_{\mathbf{Z}_{p}}\mathbf{Q}_{p})].$$
\textbf{Acknowledgement:} I would like to thank Jack Morava for his comments on the earlier version of this article and for shearing his global vision. I'm grateful to Lars Hesselholt and Bjorn Dundas for there enlightening explanations about cyclotomic trace formula. I am deeply grateful to Bertrand To\"en for his interest and remarks which were determinant for the new version of the article. 

\bibliographystyle{plain} 
\bibliography{champollion}

\end{document}